\magnification=1160

\vskip4.0truecm
\centerline {\bf THE DEGREE OF THE GENERATORS OF THE
CANONICAL RING} \vskip0.6truecm\centerline {\bf OF SURFACES OF GENERAL TYPE
WITH ${\bf p_g=0}$} \vskip0.6truecm\centerline{by}
\vskip0.6truecm\centerline{ Margarida Mendes
Lopes} \vskip1.5truecm \noindent {\bf
Introduction.} \vskip0.3truecm Upper bounds for the degree
of the generators of the canonical rings of surfaces of
general type were found by Ciliberto [C].
 In particular it was established that the canonical ring of a minimal
surface of general type with $p_g=0$ is generated by its elements of degree
lesser or equal to $6$, ([C], th. (3.6)). This was the best bound possible to
obtain at the time, since Reider's results, [R], were not yet available. In
this note, this bound is improved in
some cases (theorems (3.1), (3.2)).\par In particular it is shown  that if 
 $K^2\geq 5$, or if $K^2\geq 2$ and $|2K_S|$ is base point free this bound can
be lowered to  $4$. This result is proved by showing first that, under the
same hypothesis, the degree of the bicanonical map is lesser or equal to $4$ if $K^2\geq 3$,
(theorem (2.1)), implying that the hyperplane sections of the bicanonical
image have not arithmetic genus $0$.
 The result on the generation of the canonical ring then follows by the
techniques utilized in [C].  \vskip0.4truecm\noindent

\noindent {\bf Notation and conventions.}\vskip0.2truecm

We will denote by $S$ a projective algebraic surface over the complex field.
Usually $S$ will be smooth, minimal, of general type. \par

We denote by $K_S$, or simply
by $K$ if there is no possibility of confusion, a canonical divisor on $S$. As
usual, for any sheaf ${\cal F}$ on $S$, we denote by $h^i(S,{\cal F})$ the
dimension of the cohomology space  $H^i(S,{\cal F})$, and by $p_g$ and $q$ the
geometric genus and the irregularity of $S$.\par

By a {\it curve} on $S$ we mean an effective, non zero divisor on $S$. We will
denote the intersection number of the divisors $C$, $D$ on $S$ by $C\cdot D$
and by  $C^2$ the self-intersection of the divisor $C$. We denote by $\equiv$ the
linear equivalence for divisors on $S$. $|D|$ will be the complete linear
system of the effective divisors $D'\equiv D$, and $\phi_D:S\to {\bf
P}(H^0(S,{\cal O}_S(D)^\vee)=|D|^\vee$ the natural rational map defined by
$|D|$.\par We will denote by $\Sigma_d$ the rational ruled surface ${\bf
P}({\cal O}_{\bf P^1}\oplus {\cal O}_{\bf P^1}(d))$, for $d\geq 0$. 
$\Delta_{\infty}$ will denote the section of $\Sigma_d$ with minimum
self-intersection $-d$ and  $\Gamma$ will be a fibre of the projection to
${\bf P^1}$.  
 
\vskip0.3truecm\noindent {\it Acknowledgements.} This work was started during
my stay at the University of Tor Vergata in November 1995. I wish to thank the
University of Tor Vergata for hospitality, the Italian CNR and JNICT for
financial support.\par I wish also to thank Ciro Ciliberto and Barbara Fantechi
for useful conversations and the referee for his interest and stimulating suggestions.

  \vskip0.6truecm

\noindent {\bf 1. Auxiliary results.} \vskip0.3truecm \noindent In this section various
results that will be used in the sequel are listed.
 
\vskip0.3truecm \noindent (1.1) (Xiao [X], Reider [R])
{\sl Let $S$ be a minimal surface of general type with
$p_g=0$. If $K^2\geq 2$, $\phi_{2K}$ is generically
finite and if $K^2\geq 5$, $\phi_{2K}$ is a morphism.}

\vskip0.3truecm \noindent (1.2) (Catanese [Ca], Reider [R],
Bombieri-Catanese [BC])) {\sl Let $X$ be the canonical
model of a minimal surface $S$ of general type with
$p_g=0$. If $K_S^2\geq 3$, $\phi_{3K_X}$ is an embedding and
if $K_S^2=2$, $\phi_{3K_X}$ is a birational morphism.}  

\vskip0.3truecm \noindent (1.3) (Xiao [X], p.37) {\sl Let $S$ be a minimal surface
of general type with $p_g=0$. If $S$ contains a genus 2
pencil, then $K^2\leq 2$.}

\vskip0.3truecm \noindent (1.4) (Ciliberto, [C], th. (3.5)) {\sl The canonical
ring of a minimal surface of general type with $p_g=q=0$ is generated by its
elements of degree $\leq 6$.} 
\vskip0.3truecm \noindent (1.5) (Ciliberto, [C],
theorem (1.11), observation (1.3) ) {\sl Let $|D|$ be a complete, base point free linear series on
a smooth irreducible curve $C$. If $|D|$ is not composed with a rational
involution, or if dim$|D|=1$, then the natural map $H^0(C,{\cal O}_C(D))\otimes
H^0(C,\omega_C)\to H^0(C,{\cal O}_C(D)\otimes \omega_C)$ is surjective.}

\vskip0.3truecm \noindent (1.6) (De Franchis, [DF], see also Catanese, Ciliberto,
 [CC]) {\sl Let $X$, $Y$ be two smooth surfaces and $f:X\to Y$ be a finite
double cover. If $p_g(Y)=q(Y))=0$ and $q(Y)>0$, then the Albanese image of $X$
is a curve $C$, the involution determined by $f$ acts as (-1) on Alb(X) and
there is a morphism $h: Y\to {\bf P}^1$
 such that $f:X\to Y$ is the fiber product of $p:C\to {C\over \pm 1}\simeq
{\bf P}^1$ and $h$.}  \vskip0.3truecm \noindent (1.7)  (Arakelov, see [Be],
pg. 343) {\sl Let $S$ be a minimal surface of general type and let $f:S\to B$
be a genus $b$ pencil of curves of genus $g$. Then $K^2\geq 8(g-1)(b-1)$.}
\vskip0.3truecm \noindent (1.8) (Ciliberto, Francia, Mendes Lopes, [CFM]) Let
S be an algebraic surface and $\mu\in Pic^0(S)$. If $C$ is a curve on $S$ such
that $C^2>0$, then $\mu_{|C}$ is non-trivial.

\vskip0.5truecm \noindent {\bf 2. Degree of the bicanonical map}
 \vskip0.3truecm \noindent (2.1) {\bf Theorem.} {\sl Let $S$  be a minimal
surface of general type with $p_g=0$, $K^2\geq 3$, such
that
 $\phi_{2K}$ is a morphism. Then
degree $\phi_{2K}\leq 4$ and
 $X:=$ Im $\phi_{2K}$ is  a surface of degree bigger or equal to $n$ in ${\bf
P}^n$. }
\vskip0.2truecm\noindent {\it Proof.} By (1.1) $\phi_{2K}$ is generically
finite. Let $d$, $m$ be respectively the degrees of $\phi_{2K}$ and $X$.
Since, by assumption, $|2K|$ is basepoint free, one has $d\cdot
m=(2K)^2=4K^2$.\par Because $p_2=K^2+1$ and $X$ is not contained in any
hyperplane of ${\bf P}^{K^2}$, deg $X\geq K^2-1$. An easy calculation yields
then that $d>4$   can hold  only in the following two cases:\par\noindent (a)
$K^2=5$ and $d=5$ or \par\noindent (b) $K^2=3$ and $d=6$.\par Notice that in
these cases (and exactly in these) we have deg $X=K^2-1$. Now to prove
the theorem we are going to show that neither case can occur.\vskip0.3truecm\noindent (a) Suppose
that
$K^2=5$ and $d=5$. It is well known that a surface of degree $4$ in ${\bf
P}^5$ is one  of the following surfaces:  \par\noindent (i) $\Sigma_0$ 
embedded by
$|\Delta_{\infty}+2\Gamma|$ or $\Sigma_2$  embedded by
$|\Delta_{\infty}+3\Gamma|$;
 \par\noindent (ii)  $\Sigma_4$ mapped by $|\Delta_{\infty}+4\Gamma|$ ( i.e. a
cone over the rational curve of degree $4$);
   \par\noindent (iii) the Veronese surface. \par 

Since $X$ is the bicanonical image of $S$, $X$ must be the
Veronese surface. In fact in both cases (i) and (ii), the pull-back of
the image of $|\Gamma|$ would give a pencil $D$  on $S$ such that $2K\cdot D=5$,
an odd number, which is impossible. Hence if
degree $\phi_{2K}=5$, $X$ is the Veronese surface. 
 Denote by $L$ the image on $X$ of a line of  ${\bf P}^2$
and let $F$ be the inverse image of $L$ on $S$. Notice that, because $h^0(S,
{\cal O}_S(F)\geq 3$ and $\chi({\cal O}_S)=1$, $h^1(S,{\cal O}_S(F)\geq 2$. 
Since $2F$ is linearly equivalent to $2K$, $\eta:=K-F$ is a 2-torsion element
of $Pic(S)$. \par Let $g:\tilde S\to S$ be the \'etale double cover associated
to $\eta$. One has $K_{\tilde S}^2=10$, $\chi({\cal O}_{\tilde S})=2$
and $q(\tilde S)=q(S) + h^1(S,{\cal O}_S(F)\geq 2$. Therefore we can apply 
 (1.6) to the double cover $g:\tilde S\to S$, obtaining the existence of a
fibration $f:\tilde S\to B$, with $B$ a smooth curve of genus $b\geq 2$.
Because $K_{\tilde S}^2=10$, by (1.7) the genus of the general fibre of $f$
must be $2$. This implies in turn that also $S$ has a fibration in genus $2$
curves, contradicting (1.3). Therefore the case $K^2=5$, $d=5$ does not
occur.\vskip0.3truecm\noindent (b) If $K^2=3$ and deg $\phi_{2K}=6$, the image of $\phi_{2K}$ is
necessarily the quadric cone in ${\bf P}^3$. Otherwise the image of
$\phi_{2K}$ would be the non-singular quadric in $\bf {P}^3$ and thus
$2K\equiv D_1+D_2$, with $K\cdot D_i=3$, which is not congruent to $D_i^2=0$
(mod 2). So the image of $\phi_{2K}$ is a quadric cone and therefore we have
$2K\equiv 2D+G$, where $|D|$ is a rational pencil without fixed part,
satisfying $K\cdot D=3$ and   $G$ is an effective divisor, possibly $0$, such
that $K\cdot G=0$. Notice that, if $G\neq 0$, every irreducible component $\theta$ of $G$ is a curve such that
$\theta^2=-2$, $K\cdot \theta=0$. We can write $G$ uniquely as $G=2G'+\Gamma$,
where either $\Gamma=0$ or $\Gamma=\theta_1+...+\theta_r$, with
$\theta_1,...,\theta_r$  distinct irreducible curves. If $\Gamma\neq 0$,
following the same argument as in [CDe], prop.1.5, one can show that the
curves $\theta_1,...,\theta_r$ are disjoint, and thus the curve
$\Gamma:=\theta_1+...+\theta_r$ is smooth.\par In either case we have
$\Gamma\equiv 2\gamma$, where $\gamma\equiv K-D-G'$ and so we can consider the
smooth double cover $\pi:S'\to S$, branched over $\Gamma$ and determined by
$\gamma$. 
 By the double cover formulas one has (with r=0, if $\Gamma=0$)
\vskip0.2truecm\noindent
$\chi ( S')=2\chi(S)+{1\over 2}
\gamma\cdot (K_S+\gamma)=2-{r\over 4}$ ; 
\vskip0.2truecm\noindent$K_{S'}^2=2(K_S+\gamma)^2=2(3-{1\over 2}r)=6-r$ ; 
\vskip0.2truecm\noindent$p_g(S')=h^0(S, {\cal O_S}(K_S+\gamma))+h^0(S, {\cal
O_S}(K_S))=h^0(S,{\cal O_S} (D+G'+\Gamma)) +0\geq 2.$ \vskip0.2truecm Suppose
that
$\Gamma\neq 0$. Because
$\chi(S')\geq 1$, we must have $r=4$, yielding $\chi(S')=1$ and  $q(S')\geq
2$. In turn, by (1.6), this implies that the Albanese image of $S'$ is a curve
of genus $b=q(S')$. But then we have a contradiction to the inequality (1.7). In fact
$S'$ contains exactly four exceptional divisors, lying over the curves $Z_i$
and thus the minimal model $\tilde S$ of $S'$ has
$K_{\tilde S}^2=6< 8(q-1)$.\par Therefore $\Gamma=0$, $\gamma$ is a 2-torsion divisor
and $\pi$ is
\'etale. $S'$ is thus a minimal surface with $\chi(S')=2$ and $K_{S'}^2=6$.
Since
$p_g\geq 2$, $S'$ is irregular and in fact  $q=1$ (
$q\geq  2$  leads to the same contradiction as in the previous
paragraph). \par Let us consider now the Albanese fibration of $S'$,
$\alpha:S'\to E=Alb(S')$. By (1.6), the fixed point free involution
$i$ induced by
$\pi:S'\to S$ acts as (-1) on
$E$ and therefore there are  4 fibres of the Albanese pencil
which are stable under
$i$. Since $\pi$ is \'etale, this implies that the pencil induced on $S$ by
the Albanese pencil has 4 double fibres
   $2F_1,...,2F_4$.
The existence of these 4 double fibres implies in turn the existence in
$Pic(S)$ of seven distinct non-zero 2-torsion divisors, namely
$\eta_{ij}=F_i-F_j$, $i,j\in \{1,2,3,4\},i<j$ and $\eta=F_1+F_2-F_3-F_4$.
Notice that, given a non-trivial 2-torsion divisor $\mu $ on $S$, by the
Riemann-Roch theorem one has always $h^0(S, {\cal O}_S(K+\mu))\geq 1$.\par
Consider now
 the pencil
$|D|$. Since $K\cdot D=3$, $D^2$ is an odd number bigger than 0 and either
$G\neq 0$, $D^2=1$ and $g(D)=3$, or $G=0$, $D^2=3$ and $g(D)=4$. In either case
a general curve $D'$ in $|D|$ is smooth. Furthermore, by (1.8), the seven linear
systems
$|K+\eta_{ij}|$, $i,j\in \{1,2,3,4\},i<j$ and $|K+\eta|$ cut on $D'$ seven
distinct effective divisors $N_{ij}$ and $N$ of degree 3.\par
 Consider Im $r:H^0(S, {\cal
O}_S(2K))\to H^0(D', {\cal O}_{D'}(2K))$. Since $\phi_{2K}(D')$ is a line and
$|2K|$ is basepoint free, Im $r$  is a $g^1_6$ without base points. Now $2N_{ij}$, for $i,j \in \{1,2,3,4\}, i<j$ and
$2N$ belong to Im $r$ and each of these divisors gives a contribution of at
least $3$ for the degree of the ramification divisor $R$ of the 6-1 morphism
$D'\to {\bf P}^1$. Therefore deg $R\geq 21$, and so by the Hurwitz formula one
has $2g(D')-2\geq 6(-2)+21=9$, i.e. $g(D')\geq 6$, which contradicts
$g(D')\leq 4$. Thus also the case $K^2=3$, $d=6$ cannot occur. $\diamondsuit$ 

\vskip0.4truecm \noindent (2.2) {\bf Observations.} \par\noindent 1. If
$K^2=2$ and $|2K|$ is base-point free, then the degree of the bicanonical map
is 8.
 \par\noindent 2. The degree of the bicanonical map of the Burniat surface
with $K^2=3$ (cf. [P]) is 4 and therefore, at least for $K^2=3$, the bound
established in theorem (2.1) is sharp.

\vskip0.5truecm\noindent {\bf 3. The degree of the generators of the canonical
ring}
\vskip0.3truecm   
 Let us fix some notation first. $C$ will be a general element in $|2K_S|$, whilst $H_m:=H^0(C,
{\cal O}_C(mK_S))$ and  $R_m:= H^0(S, {\cal O}_S(mK_S))$. We will  denote
by $R$ the canonical ring of $S$ and by $H$ the ring $H:=\oplus_{m=1}^{\infty}
H_m$. \par   Notice that,
because $q=0$, the restriction maps $r_m:R_m\to H_m$ are surjective for every
$m\in {\bf N}$, and so the natural graded morphism $r:R\to H$ induced by the
restriction maps is also surjective. Furthermore ker $r$ is a principal ideal
of $R$, generated by one element of degree 2 (cf.[C], \S 2). By (1.4), $R$ is generated by its
elements of degree lesser or equal to $6$. To lower this bound to $5$ it will
be enough to show that the map $H_3\otimes H_3\to H_6$ is surjective, and
similarly to lower the bound to $4$ it will be enough to show that  the map
$H_2\otimes H_3\to H_5$ is surjective (cf.[C], \S 2).  

\vskip0.3truecm \noindent (3.1) {\bf Proposition.} {\sl Let $S$ be a minimal
surface of general type with $p_g=0$, $K^2\geq 2$ such that $|2K|$ has no
fixed components. Then the canonical ring of $S$ is generated
by its elements of degree $\leq 5$.} 

\vskip0.2truecm \noindent {\it Proof.} By (1.1)
Im$\phi_{2K}$ is a surface and so by our hypothesis the
general curve $C\in |2K|$ is irreducible. Since, by (1.2), 
$\phi_{3K}$ is birational, $C$ is not
hyperelliptic and therefore the map
$$H^0(C, \omega_C)\otimes H^0(C, \omega_C)\to H^0(C,
\omega_C^{\otimes 2})$$ is surjective. Since, by adjunction, $\omega_C\simeq
{\cal O}_C(3K_S)$, the map $H_3\otimes H_3\to H_6$ is surjective. Hence, by
(1.4) and the previous considerations we have the assertion.$\diamondsuit$

\vskip0.3truecm \noindent (3.2) {\bf Proposition.} {\sl Let $S$ be a minimal
surface of general type with $p_g=q=0$. The canonical ring of $S$ is generated
by its elements of degree $\leq 4$ if $S$ satisfies one of the following conditions:

\vskip0.1truecm\noindent i) $K^2\geq 5$; or  \vskip0.1truecm\noindent ii) 
 $|2K|$ is base point free and $K^2\geq 2$. } \vskip0.2truecm \noindent {\it
Proof.} By (1.4)and proposition (3.1), it suffices to show that the map
$H_2\otimes H_3\to H_5$ is surjective. Let
$C$ be a general bicanonical curve. Since $|2K|$ is base point free (in case
(i) by (1.1)), the linear series $|{\cal O}_C(2K_S)|$ is base point free. \par
If $K^2\geq 3$, by  theorem (2.1), Im$\phi_{2K}$ is not a surface of
minimal degree $n-1$ in ${\bf P}^n$. Therefore $|{\cal O}_C(2K_S)|$ is not composed
with an rational involution and so, by (1.5), we have the surjectivity of the
above map.\par If $K^2=2$, one has $h^0(S, {\cal O}_S(2K))=3$ and thus
dim$|{\cal O}_C(2K_S)|=1$, implying, again by (1.5), that the map $H_2\otimes
H_3\to H_5$ is surjective.$\diamondsuit$   \vskip0.3truecm \noindent (3.3)
{\bf Observation.}
\par\noindent  To my knowledge, all the known examples of minimal surfaces
of general type with $p_g=q=0$ and $K^2=2,3,4$ have bicanonical system free
from base points.

\vskip0.5truecm \noindent {\bf References.}

\vskip0.3truecm

\item{[BC]} E. Bombieri and F. Catanese, {\it The
tricanonical map of a surface with $K^2=2$, $p_g=0$}, C.P.
Ramanujam - A tribute, Bombay, 1978.

\item{[Be]} A. Beauville, {\it L'in\'egalit\'e $p_g\geq 2q-4$ pour les surfaces de
type g\'en\'erale}, Bull. Soc. Math. de France, 110 (1982), 343-346 
\item{[C]} C. Ciliberto, {\it Sul grado dei generatori dell'anello canonico di una
superficie di tipo generale}, Rend. Sem. Mat. Univers. Politecn. Torino, 41, 3 (1983), 83-111 

\item {[Ca]} F. Catanese, {\it
Footnotes to a theorem of I. Reider}, Algebraic Geometry, Lecture Notes in
Math, 1417, Springer,1990.

\item {[CC]} F. Catanese, C. Ciliberto, {\it On the irregularity of cyclic
coverings of algebraic surfaces},
 in Proceedings of the conference "Geometry of complex projective varieties", Cetraro, 1990, Mediterranean Press.
\item {[CDe]} F. Catanese, O.Debarre, {\it Surfaces with $K^2=2$, $p_g=1$, $q=0$}, J.Reine angew. Math, 395 (1989), 1-55

\item {[CFM]} C. Ciliberto, P.Francia, M. Mendes Lopes, {\it Remarks on the bicanonical map for surfaces of general type}, to be published in Math.Zeit.
\item {[DF]} M. De Franchis, {\it Sugl'integrali di Picard
relativi a una superficie doppia}, Rend. Circ. Mat.
Palermo XX (1905), 331-334
\item {[P]} C.A.M. Peters, {\it On certain examples of surfaces with $p_g=0$ due
to Burniat}, Nagoya Math. J., 66 (1977), 109-119 
\item {[R]} I. Reider, {\it
Vector bundles of rank 2 and linear systems on algebraic surfaces}, Ann. of
Math. 127 (1988), 309-316 

\item {[X]} G. Xiao, {\it
Finitude de l'application canonique des surfaces de type
g\'en\'eral}, Boll. Soc. Math. France, 113 (1985), 23-51

\vskip0.5truecm\noindent Author's adress:\hfill e-mail\par\noindent
C.A.U.L.\hfill mmlopes@ptmat.lmc.fc.ul.pt\par\noindent Av.Prof.Gama Pinto,2
\par\noindent 1699 Lisboa Codex
\par\noindent Portugal
\bye